\documentclass{amsart}
\usepackage{amsmath}
\usepackage{amssymb}
\usepackage{todonotes}
\usepackage{tikz}
\usepackage{pgfplots}
\pgfplotsset{compat=1.18}
\usepackage{epsfig}
\usepackage{tabularray}
\usepackage{url}
\usepackage{hyperref}
\usepackage{cleveref}
\usepackage{array}
\usepackage{subcaption}
\newtheorem{theorem}{Theorem}
\theoremstyle{remark}
\newtheorem{remark}[theorem]{Remark}
\newtheorem{definition}[theorem]{Definition}
\numberwithin{theorem}{section}

\title{Optimal Play, Nontransitivity, and Nash Equilibria in Dice Bingo}
\date{\today}

\author{David J.~Hemmer}
\address{Department of Mathematical Sciences\\
Michigan Technological University\\
Houghton, MI 49931}
\email{djhemmer@mtu.edu}

\author{Benjamin W.~Ong}
\address{Department of Mathematical Sciences\\
Michigan Technological University\\
Houghton, MI 49931}
\email{ongbw@mtu.edu}

\newlength{\bingosize}
\newcommand{\bingo}[9]{%
\begin{tikzpicture}[x=\bingosize,y=\bingosize]
  \draw[thick] (0,0) grid (3,3);

  \node at (0.5,2.5) {#1};
  \node at (1.5,2.5) {#2};
  \node at (2.5,2.5) {#3};

  \node at (0.5,1.5) {#4};
  \node at (1.5,1.5) {#5};
  \node at (2.5,1.5) {#6};

  \node at (0.5,0.5) {#7};
  \node at (1.5,0.5) {#8};
  \node at (2.5,0.5) {#9};
\end{tikzpicture}
}

\begin{document}

\maketitle

\begin{abstract}
  
We study Dice Bingo, a game in which players fill a $3\times3$ bingo
board whose entries are possible sums of two fair dice. After each
roll, a player marks one matching square, and the goal is to complete
a row, column, or diagonal. We model optimal play for a fixed board as
a finite Markov decision process and derive Bellman equations that
compute the exact expected number of rolls required to obtain a
bingo. Using this framework, we identify a unique optimal board up to
natural symmetries and determine its exact expected completion time.

We then investigate head-to-head competition in which two players
observe the same sequence of dice rolls. By analyzing a joint Markov
chain that tracks both boards simultaneously, we compute (in exact
arithmetic) win, loss, and tie probabilities. Surprisingly, a board
with a worse expected completion time can nevertheless be favored in
head-to-head competition. Motivated by this phenomenon, we exhibit
nontransitive triples of bingo boards: board $A$ is favored against
board $B$, board $B$ is favored against board $C$, and board $C$ is
favored against board $A$.

Finally, we consider strategic play in which players adapt their
choices to their opponent's board rather than merely minimizing their
own completion time. In this setting, optimal decisions depend on the
opponent's state, leading naturally to game-theoretic analysis. We
present a position with no pure Nash equilibrium and compute an
explicit mixed Nash equilibrium.
\end{abstract}

\vspace{10px}

\noindent{\textbf{Keywords: }{Dice Bingo, Markov Decision Process, Bellman Equation, Nash Equilibrium}} 

\vspace{10px}

\section{Introduction}

The authors' university sponsors a monthly math circle for students in
grades 2-5. A recent activity included a game from the Julia Robinson
Math Festival \cite{jrmf_dice_bingo} called ``Dice Bingo''. Students
discover that when two dice are rolled, the sums vary between 2 and
12. The students design individual $3 \times 3$ boards whose entries
are chosen from $\{2,3,\ldots, 12\}$, with repetitions allowed. Two
dice are repeatedly rolled, and whenever the sum matches an unmarked
entry, the student may cover one square. The goal is to win \`{a} la
tic-tac-toe, i.e., to obtain three filled squares in a row, column or
diagonal, which we call in this paper a ``bingo.''

One precocious student announced that she knew seven was the most
likely sum with two dice, and proudly filled her board in with nine
sevens! Watching this poor student lose over and over again to
students who had seemingly filled in their boards at random is what
spurred this investigation into dice bingo strategy. For the $3 \times
3$ case with two dice, we find through an exhaustive computational
search that there is a unique board (up to symmetries induced by the
dihedral group and by swapping numbers $i$ with the equally likely sum
$14-i$) with the minimal expected number of rolls to win. Surprisingly
to us, this board does not place seven in the middle square. We also
find that a magic square with sum equal to 21 is nearly optimal,
despite containing the relatively unlikely sums 3 and 11.

Inspired by Efron's nontransitive dice (popularized in a column of
Martin Gardner \cite{Gardner1970NontransitiveDice}), we also find a
set of ``nontransitive boards'' $A,B,C$.  In head-to-head competition
where each player uses solo optimal strategy to minimize their winning
times (ignoring their opponent's board), player $A$ beats player $B$,
player $B$ beats $C$, and player $C$ beats player $A$.

Lastly, we also consider a head-to-head competition in which the
players adjust their strategies based on their opponent's boards and
choices, rather than simply optimizing their own winning time. Because
each player's best response depends on the opponent's board and
choices, deterministic strategies can be exploited, leading naturally
to a mixed-strategy Nash equilibrium.

We believe dice bingo is an excellent problem to use pedagogically, as
a thorough treatment involves concepts from linear algebra,
probabilities, finite Markov chains, Bellman equations, group actions
and game theory, as well as tractable and interesting computational
problems. We also encounter open problems that are beyond our
computational capability.

\section{Expected number of rolls to completion for a single board}
\label{sec: singleboard}

We wish to determine the expected number of rolls for a single board
to win, i.e., complete a bingo. A $3\times 3$ bingo board,
\cref{fig:bingo_card},
\begin{figure}[htbp]
  \centering
  \bingo{$b_0$}
        {$b_1$}
        {$b_2$}
        {$b_3$}
        {$b_4$}
        {$b_5$}
        {$b_6$}
        {$b_7$}
        {$b_8$}
  \caption{A $3\times 3$ bingo card has entries $b_i\in \{2, 3,
    \ldots, 12\}$. We label the cells from left to right, top to
    bottom, starting the enumeration at 0.}
  \label{fig:bingo_card}
\end{figure}
is represented by the vector
\hbox{$B = (b_0, b_1, \ldots, b_8)$} where each $b_i \in \{2, 3, \ldots, 12\}$
is a possible dice sum written in cell~$i$.
Define a \emph{state} to be a subset $S \subseteq \{0,1,\ldots,8\}$
recording which cells have been marked. There are $2^9 = 512$ possible
states, and a game starts with the board in state $S=\varnothing$. The
eight winning lines, corresponding to cells in the three rows, three
columns, and two main diagonals, are
\begin{align*}
  \mathcal{L} = \bigl\{
  \{0,1,2\},\;\{3,4,5\},\;\{6,7,8\},\;
  \{0,3,6\},\;\{1,4,7\},\;\{2,5,8\},\;
  \{0,4,8\},\;\{2,4,6\}
  \bigr\}.
\end{align*}
A state $S$ is \emph{winning} if $\ell \subseteq S$ for some $\ell \in
\mathcal{L}$.  Let $\mathcal{W}$ denote the set of absorbing (winning)
states and let $\mathcal{T}$ denote the set of \emph{transient}
(non-winning) states.  One can verify by direct enumeration or using
the inclusion-exclusion principle that $|\mathcal{W}| = 282$ and
$|\mathcal{T}| = 230$.

\Cref{tbl:sum_prob} shows $p(v)$, the probability of rolling sum $v$
with two fair dice.
\begin{table}[h!]
  \centering
  \begin{tabular}{c|ccccccccccc}
    $v$ & 2 & 3 & 4 & 5 & 6 & 7 & 8 & 9 & 10 & 11 & 12 \\
    \hline
    $36 \cdot p(v)$ & 1 & 2 & 3 & 4 & 5 & 6 & 5 & 4 & 3 & 2 & 1
  \end{tabular}
  \caption{Probabilities of rolling sum $v$ with two fair dice.}
  \label{tbl:sum_prob}
\end{table}

\subsection{The Markov decision process}

On each roll, the player observes the dice sum~$v$ and may mark one
unmarked cell whose value equals~$v$ (if such a cell exists). When
multiple unmarked cells share the value~$v$, the player \emph{chooses}
which cell to mark. The goal is to reach a winning state in as few
rolls as possible. The game can be viewed as a finite-state decision
process. At each state the dice introduce randomness, while repeated
values on the board may require the player to choose among several
possible actions. This naturally leads to a Markov decision process
(MDP).

The MDP is monotone and acyclic up to self-loops, as every non-self
transition strictly increases the number of marked cells. This acyclic
structure allows us to solve the MDP using \emph{backward induction},
where states are processed in decreasing order of $|S|$, starting from
the winning states.

\subsection{Bellman equation}

Since the outcome of each roll depends only on the current set of
marked squares, and since the player may sometimes choose among
several possible marks, the optimal expected time to completion is
described by Bellman equations for the corresponding Markov decision
process. Write $V(S)$ for the expected number of rolls from state~$S$
to a winning state under optimal play,
\begin{align*}
  V(S) = \mathbb{E}[\text{rolls remaining to win | current state } S]
\end{align*}
For a winning state $S \in \mathcal{W}$, we have $V(S) = 0$ since no
further rolls are needed (a bingo has already been found). For a
transient state $S \in \mathcal{T}$, the expected number of rolls
under optimal play satisfies a Bellman equation, which we now derive.

Given state $S$, define the set of \emph{hitting sums}:
\begin{align*}
  H(S) = \{v \in \{2,\ldots,12\} : \text{there exists } i \notin S
  \text{ with } b_i = v\}.
\end{align*}
These are the dice sums that match at least one unmarked cell.  Any
sum $v \notin H(S)$ is \emph{wasted}: the state remains at~$S$ if a
$v$ is rolled. The probability of a useful roll is
\begin{align*}
  p_{\mathrm{hit}}(S) = \sum_{v \in H(S)} p(v).
\end{align*}
For each hitting sum $v \in H(S)$, the player chooses which cell to
mark among the candidates
\begin{align*}
  C(S, v) = \{i \notin S : b_i = v\}.
\end{align*}
Marking cell $c \in C(S,v)$ transitions to state $S \cup \{c\}$.  When
there is a choice, the player will mark the cell that minimizes the
expected number of future rolls:
\begin{align*}
  c^*(S,v) = \operatorname*{arg\,min}_{c \,\in\, C(S,v)} V(S \cup \{c\}),
\end{align*}
where $V(S \cup \{c\}) = 0$ if $S \cup \{c\} \in \mathcal{W}$.  From
state $S$ we have:
\begin{itemize}
\item With probability $p(v)$ for each $v \in H(S)$, the roll is
  useful and we transition to $S \cup \{c^*(S,v)\}$.
\item With probability $1 - p_{\mathrm{hit}}(S)$, the roll is wasted
  and we remain at $S$.
\end{itemize}
We obtain:
\begin{equation}
  \label{eq:bellman_raw}
  V(S) = \underbrace{1}_{\text{current roll}} + \underbrace{\sum_{v
      \in H(S)} p(v)\, V\!\bigl(S \cup
    \{c^*(S,v)\}\bigr)}_{\text{useful roll}} + \underbrace{\bigl(1 -
    p_{\mathrm{hit}}(S)\bigr) \, V(S)}_{\text{wasted roll}}.
\end{equation}
The ``$1$'' in \cref{eq:bellman_raw} accounts for the current roll and
the next two terms account for the expectations after a hitting sum or
a wasted roll respectively.  Rearranging,
\begin{align*}
  V(S) - \bigl(1 - p_{\mathrm{hit}}(S)\bigr)\, V(S) = 1 + \sum_{v \in
    H(S)} p(v)\, V\!\bigl(S \cup \{c^*(S,v)\}\bigr),
\end{align*}
and simplifying,
\begin{align*}
  p_{\mathrm{hit}}(S)\, V(S) = 1 + \sum_{v \in H(S)} p(v)\, V\!\bigl(S
  \cup \{c^*(S,v)\}\bigr),
\end{align*}
gives
\begin{equation}
  \label{eq:bellman}
  V(S) = \frac{1 + \displaystyle\sum_{v \in H(S)} p(v)\, V\!\bigl(S
    \cup \{c^*(S,v)\}\bigr)}{p_{\mathrm{hit}}(S)}.
\end{equation}

Notice that \cref{eq:bellman} expresses the expected number of rolls
$V(S)$ to win from state $S$ in terms of various $V(S')$ where $S'$
are states with strictly more squares marked. This lets us recursively
obtain a closed form for $V(S)$ starting with $V(S)=0$ for winning
states and working backwards. We make this precise in
\Cref{sec:backward_induction}.

\subsection{Backward induction}
\label{sec:backward_induction}

Since marking a cell increases $|S|$ by exactly one, the value $V(S
\cup \{c\})$ on the right-hand side of \cref{eq:bellman} involves a
state with $|S|+1$ marked cells.  Therefore, if we process states in
decreasing order of $|S|$, every value on the right-hand side is
already computed when we evaluate $V(S)$. So we can compute any $V(S)$
via the following recursive procedure, noting that if $|S| \geq 7$ it
must already be winning:

\begin{enumerate}
\item \textbf{Base case.} For all $S \in \mathcal{W}$, set $V(S) = 0$.

\item \textbf{Recursive step.} For $k = 6, 5, \ldots, 1, 0$, process
  each transient state $S$ with $|S| = k$:

  \begin{enumerate}
  \item Compute $H(S)$ and $p_{\mathrm{hit}}(S)$.
  \item For each $v \in H(S)$, find $c^*(S,v)$ by evaluating $V(S \cup
    \{c\})$ for all $c \in C(S,v)$ and choosing the minimum.  (All
    these values are already known because $|S \cup \{c\}| = k + 1 >
    k$.)
  \item Compute $V(S)$ via \eqref{eq:bellman}.
  \end{enumerate}

\item \textbf{Output.} The quantity of interest is $V(\varnothing)$:
  the expected number of rolls from the initially unmarked board.
\end{enumerate}

\begin{remark}
  All arithmetic is over $\mathbb{Q}$ so $V(\varnothing) \in
  \mathbb{Q}$, and the computation yields an \emph{exact} rational
  answer.
\end{remark}

Once the backward induction is complete, the optimal strategy for the
given board is recorded as a partial function $c^*\colon \mathcal{T}
\times \{2,\ldots,12\} \to \{0,\ldots,8\}$ specifying, for each state
and each hitting sum, which cell to fill in. This fully determines the
player's optimal strategy. Note we also obtain, during the course of
the computation, the expected winning time and optimal strategy
starting with any partially marked board. Also note in step 2(b) above
there may be ties, in which case we could simply specify to mark the
lowest numbered cell from among the equal choices. This observation
will be important later when we consider head to head matchups, where
the choice may matter.

\begin{remark}
  The optimal strategy $c^*$ can be unintuitive at first glance if
  the Bellman equation is not solved. For example, consider the bingo
  card below, which will be utilized later in
  \cref{sec:slower}. Suppose that we have the following  unmarked board,
  \begin{center}
    \bingo{6}{7}{6}{7}{7}{7}{6}{6}{6} 
  \end{center}
  and a hitting sum $v=6$
  is rolled. Which cell is marked in an optimal strategy? We leave it
  as an exercise to the reader to show that cell 7 is the optimal
  square to mark, counter to intuition.
\end{remark}

\begin{remark}
  In lieu of solving the Bellman equation for the optimal strategy
  given a choice of unmarked squares, one could instead choose
  randomly according to a fixed probability distribution which square to fill, or apply a fixed
  strategy. For example, preferentially fill cell 4 (middle), then
  corners (cells 0,2,6,8), then remaining cells.
\end{remark}

\subsection{Board of all sevens}

Consider the board $B = (7,7,7,7,7,7,7,7,7)$ that our unlucky student
constructed. Every cell has value~$7$, so $H(S) = \{7\}$ for all $S
\in \mathcal{T}$, and $p_{\mathrm{hit}}(S) = p(7) = 1/6$.

Any three cells in a line suffice.  The earliest bingo occurs after
marking exactly three cells that form a winning line.  The optimal
strategy must mark cells to complete a line in exactly three useful
rolls. Since each useful roll arrives after a geometric wait with
parameter $1/6$, the expected number of rolls is
\begin{align*}
  V(\varnothing) = \underbrace{6}_{\text{wait for 1st seven}} +
  \underbrace{6}_{\text{wait for 2nd seven}} +
  \underbrace{6}_{\text{wait for 3rd seven}} = 18.
\end{align*}
As we will see the best boards have expected times substantially less
than 18, which explains our student's difficulties winning.

\subsection{Optimal board}

\Cref{fig:optimal} displays one of the 64 equivalent (up to symmetry)
boards that we found with the smallest expected win time. Notice we
can apply any of the four rotations and four reflections in the
dihedral group $D_4$ to obtain boards with the same expected win
time. We can also swap 6's and 8's, 5's and 9's or 4's and 10's which
gives $8\cdot 2^3=64$ equivalent boards.
\begin{figure}[htbp]
  \centering
  \bingo{8}{8}{9}{7}{6}{10}{7}{4}{5}
  \caption{Optimal board with $V(\varnothing) \approx
    6.3094$. Surprising to us, $b_4 \neq 7$; intuition expects that
    dice sum with the highest probability would be in the center cell
    since it affects the most number of winning lines. However, this
    is not the case since it is non-optimal for too many lines to
    depend on a seven.}
  \label{fig:optimal}
\end{figure}
The optimal board\footnote{The optimal board was identified by
computing the exact expected number of rolls for all $11^9 =
2,357,947,691$ possible boards using over 2000 CPU hours on our local
HPC cluster. The code does not exploit symmetry. The python script,
which can be compiled for optimal runtime, is provided for
reproducibility. See \url{https://github.com/ongbw/dice-bingo}}  contains seven distinct values
$\{4,5,6,7,8,9,10\}$, so
\begin{align*}
  p_{\mathrm{hit}}(\varnothing) &= p(4)+p(5)+p(6)+p(7)+p(8)+p(9)+p(10)
  \\ &= \frac{3+4+5+6+5+4+3}{36} = \frac{30}{36} = \frac{5}{6}.
\end{align*}
That is, $5/6$ of all rolls are useful from the unmarked board: only sums
of $2$, $3$, $11$, or $12$ are wasted.  Backward induction through the
$230$ transient states yields the exact expected number of rolls:
\begin{align*}
  V_{\text{optimal}}(\varnothing) =
  \frac{47{,}546{,}657{,}067{,}260{,}786{,}722{,}139}{7{,}535{,}828{,}431{,}282{,}951{,}800{,}000}
  \approx 6.3094.
\end{align*}
This is the minimum over all $11^9 = 2{,}357{,}947{,}691$ possible
boards.  The optimal strategy $c^*$ determines, for each of the $230$
transient states and each useful dice sum, which cell to fill in.
\begin{remark}
  Although our brute-force search for an optimal $3\times 3$ board was
  successful, this approach becomes quickly intractable even for the
  $4 \times 4$ board, where there are now $11^{16}$ possible
  boards. Even with careful optimization for symmetry and heuristics
  to eliminate obviously poor choices, identifying the $4 \times 4$
  board using a brute force approach is an impractical use of
  computational resources.
 \end{remark}

\section{Connection to absorbing Markov chains}
\label{sec:MarkovChain}

The Bellman computation gives the optimal strategy; once that strategy
is fixed, the MDP becomes an ordinary absorbing Markov chain.  Suppose
we fix a starting board. We first determine the optimal strategy $c^*$
via backward induction, as in the previous section. Now for each state
and hitting roll we know the probability of that roll and which state
we transition to if that roll arises, our ``strategy" is completely
specified. That is, for each board we simply have an ordinary Markov
chain with 512 states and transition matrix
\begin{align*}
  T = \begin{bmatrix} Q & R \\ 0 & I \end{bmatrix},
\end{align*}
where $Q$ is the $|\mathcal{T}| \times |\mathcal{T}|$ matrix of
transitions among transient states, $R$ is the $|\mathcal{T}| \times
|\mathcal{W}|$ matrix of absorption transitions, and $I$ is the
identity on the absorbing (winning) states. As we mentioned earlier,
the underlying graph for the matrix $Q$ is directed and acyclic except
for possible loops. The $(i,j)$ entry in $Q^k$ gives the probability
that starting in transient state $i$ we arrive in state $j$ after $k$
steps. Thus the matrix $$I + Q + Q^2 + Q^3 + \cdots = (I-Q)^{-1}$$
keeps track of the total number of times the chain is expected to
visit a transient state $j$ before absorption. Specifically, the
\emph{fundamental matrix} of the Markov chain is defined to be $N = (I
- Q)^{-1}$ and entry $N_{ij}$ gives the expected number of times the
chain visits transient state~$j$ before absorption, starting from
transient state~$i$.  Summing over the transient states we get that
the expected number of rolls to absorption from state~$i$ is
\begin{equation}
  \label{eq:expectedviaMarkov}
  V(i) = \sum_{j \in \mathcal{T}} N_{ij} = \bigl[(I-Q)^{-1}\, \mathbf{1}\bigr]_i,
\end{equation}
where $\mathbf{1}$ is the all-ones vector.  This is exactly the vector
of values $V(S)$ computed by backward induction. Setting $i=0$ in
\eqref{eq:expectedviaMarkov} specifies the empty state (unmarked
board) and gives the expected number of rolls for the board to win.

The Bellman recurrence, \cref{eq:bellman}, is equivalent to the linear
system $(I - Q)\,\mathbf{t} = \mathbf{1}$, and backward induction
solves it by back-substitution, exploiting the fact that $Q$ is upper
triangular when states are ordered by decreasing number of marked
squares.

In \Cref{sec:h2h} when we look at head to head matches between two
boards, this Markov chain interpretation will be useful.

\section{Head-to-head matchups}
\label{sec:h2h}

Suppose two players, each with their own board, observe the same
sequence of dice rolls. Each player independently follows their own
optimal solo strategy $c^*$ (as computed by backward induction in
\Cref{sec: singleboard}).  Players do not see each other's boards and
do not adjust their strategies based on their opponent's state. The
game ends when at least one player completes a bingo.

Define the random variables
\begin{align*}
  T_A = \text{first roll that player A wins} \\
  T_B = \text{first roll that player B wins.}   
\end{align*}    
Because both players see the same sequence of dice rolls, $T_A$ and
$T_B$ are \emph{not} independent -- their joint distribution is
determined by the shared dice sequence.  The three outcomes of
interest are
\begin{align*}
  P(T_A < T_B), \qquad P(T_B < T_A), \qquad P(T_A = T_B).
\end{align*}
We say \textbf{$A$ beats $B$} if $P(T_A < T_B) > P(T_B < T_A)$.  These
probabilities cannot be determined only from $V_A(\varnothing)$ and
$V_B(\varnothing)$, the expected number of dice rolls to completion from the
respective unmarked boards; we need the full joint distribution, which
we now compute in exact arithmetic using a Markov chain where the
states keep track of both boards simultaneously.

\subsection{Joint state space}
The joint state of the game is the pair $(S_A, S_B)$, where $S_A$ and
$S_B$ are the states corresponding to the marked cells on boards $A$
and $B$ respectively.  A joint state is absorbing if at least one
board is in a winning state.  We classify absorbing states into three
types:
\begin{itemize}
\item \textbf{$A$ wins:} $S_A \in \mathcal{W}$ and $S_B \notin
  \mathcal{W}$ (only $A$ has completed a bingo).
\item \textbf{$B$ wins:} $S_B \in \mathcal{W}$ and $S_A \notin
  \mathcal{W}$ (only $B$ has completed a bingo).
\item \textbf{Tie:} $S_A \in \mathcal{W}$ and $S_B \in \mathcal{W}$
  (both complete a bingo on the same roll).
\end{itemize}

The transient joint states are the $230^2=52{,}900$ pairs $(S_A, S_B)
\in \mathcal{T} \times \mathcal{T}$.  For each absorbing state we
solve the linear system:
\begin{equation}
  \label{eq:h2h}
  (I - Q)\,\mathbf{x} = \mathbf{b},
\end{equation}
where $Q$ is the $52{,}900 \times 52{,}900$ transition matrix among
transient joint states and $\mathbf{b}$ is the vector of one-step
absorption probabilities into that absorbing type.  The solution at
$(\varnothing, \varnothing)$ gives the desired probabilities.  In
practice, most of these $52{,}900$ states are unreachable from
$(\varnothing, \varnothing)$ and can never arise in an actual game; an
initial breadth-first search can identify the reachable states
resulting in a much smaller linear system to solve.

\subsection{Joint transitions}

Since both boards observe the same rolls, each dice sum $v$ (with
probability $p(v)$) updates both boards simultaneously.  Each player
marks a cell according to their own pre-computed optimal strategy
playing alone, so the transition from $(S_A, S_B)$ on a roll of~$v$ is
deterministic:
\begin{align*}
  (S_A, S_B) \;\xrightarrow{\;v\;}\; \bigl(\mathrm{next}_A(S_A, v),\;
  \mathrm{next}_B(S_B, v)\bigr),
\end{align*}
where $\mathrm{next}_A(S_A, v) = S_A \cup \{c^*_A(S_A, v)\}$ if $v \in
H_A(S_A)$, and $\mathrm{next}_A(S_A, v) = S_A$ otherwise (and
similarly for board $B$).  The resulting joint state is then either
absorbing (at least one board wins) or another transient state.

\subsection{Absorption probabilities via linear system}

We enumerate the reachable transient joint states by breadth-first
search from $(\varnothing, \varnothing)$.  For each transient state
$(S_A, S_B)$, we record:
\begin{itemize}
\item the transition probabilities to other transient states, forming
  the sub-stochastic matrix $Q$;
\item the one-step absorption probabilities $b_A(S_A, S_B)$, $b_B(S_A,
  S_B)$, and $b_{\mathrm{tie}}(S_A, S_B)$ into each of the three
  absorbing types.
\end{itemize}

For each terminal type, we solve the linear system
\begin{equation}
  \label{eq:h2hrepeat}
  (I - Q)\,\mathbf{x} = \mathbf{b},
\end{equation}
where $\mathbf{b}$ is the vector of one-step absorption probabilities
into that type.  The solution $\mathbf{x}$ gives the total absorption
probability from each transient state; evaluating at the starting
state $(\varnothing, \varnothing)$ yields the exact win, loss, and tie
probabilities for the pair of boards $(A,B)$.

\begin{remark}
  As in the single-board case, one could alternatively solve
  \eqref{eq:h2hrepeat} by backward induction, processing joint states
  in decreasing order of the total number of marked boxes $|S_A| +
  |S_B|$.  Since every non-self-loop transition increases $|S_A| +
  |S_B|$ by either one or two, this ordering ensures all successor
  values are known. We use the linear system formulation here because
  it is slightly simpler to implement for the joint chain.
\end{remark}

\begin{remark}
  All arithmetic is again over $\mathbb{Q}$, so the win/loss/tie
  probabilities are rational numbers summing to~$1$.
\end{remark}

\subsection{A surprising result: ``slower'' boards may be better}
\label{sec:slower}

\begin{figure}[htbp]
\centering
  \begin{subfigure}[b]{0.48\textwidth}
  \centering
  \begin{tikzpicture}
    \begin{axis}[
      width=\linewidth,
      title={PDF},
      xlabel={Roll number ($k$)},
      ylabel={$P(T = k)$},
      xmin=1, xmax=25,
      ymin=0, ymax=0.10,
      ytick={0,0.02,0.04,0.06,0.08,0.10},
      yticklabels={0,0.02,0.04,0.06,0.08,0.10},scaled y ticks=false,
      grid=both,
      major grid style={dotted, gray!50, line width=0.5pt},
      minor grid style={dotted, gray!20},
      minor tick num=1,
      legend style={at={(1,0.97)}, anchor=north   east, font=\small},
    ]

    \addplot[blue, thick, mark=*, mark size=1.5pt] coordinates {
      (3,0.01680384) 
      (4,0.03709419) 
      (5,0.05448229) 
      (6,0.06664291)
      (7,0.07339827) 
      (8,0.07554617) 
      (9,0.07420017) 
      (10,0.07045335)
      (11,0.06523252) 
      (12,0.05925645) 
      (13,0.05304558) 
      (14,0.04695305)
      (15,0.04120110) 
      (16,0.03591503) 
      (17,0.03115174) 
      (18,0.02692214)
      (19,0.02320790) 
      (20,0.01997357) 
      (21,0.01717493) 
      (22,0.01476463)
      (23,0.01269579) 
      (24,0.01092410) 
      (25,0.00940898)
    };
    \addlegendentry{Board $A$}
    \addplot[red, thick, mark=square*, mark size=1.5pt] coordinates {
      (3,0.01502486) 
      (4,0.03542774) 
      (5,0.05493567) 
      (6,0.07021545)
      (7,0.08004994) 
      (8,0.08454327) 
      (9,0.08450019) 
      (10,0.08099962)
      (11,0.07513251) 
      (12,0.06786121) 
      (13,0.05996045) 
      (14,0.05200812)
      (15,0.04440323) 
      (16,0.03739604) 
      (17,0.03112108) 
      (18,0.02562810)
      (19,0.02090841) 
      (20,0.01691595) 
      (21,0.01358331) 
      (22,0.01083325)
      (23,0.00858666) 
      (24,0.00676760) 
      (25,0.00530632)
    };
    \addlegendentry{Board $B$}
    \end{axis}
  \end{tikzpicture}
  \caption{Solo play: independent Markov chain}
  \label{fig:solo}
  \end{subfigure}
  \hfill
  \begin{subfigure}[b]{0.48\textwidth}
    \centering
    \begin{tikzpicture}
      \begin{axis}[
          width=\linewidth,
          title={CDF},
          xlabel={Roll number ($k$)},
          ylabel={$P(T \le k)$},
          xmin=1, xmax=25,
          ytick={0,0.25,0.5,0.75,1},
          grid=both,
          major grid style={dotted, gray!50, line width=0.5pt},
          minor grid style={dotted, gray!20},
          minor tick num=1,
          legend style={at={(0.97,0.37)}, anchor=north   east, font=\small},
        ]
        
        \addplot[blue, thick, mark=*, mark size=1.5pt] coordinates {
          (3,0.01680384)
          (4,0.05389803)
          (5,0.10838032)
          (6,0.17502323)
          (7,0.24842150)
          (8,0.32396767)
          (9,0.39816784)
          (10,0.46862119)
          (11,0.53385371)
          (12,0.59311016)
          (13,0.64615574)
          (14,0.69310879)
          (15,0.73430989)
          (16,0.77022492)
          (17,0.80137666)
          (18,0.82829880)
          (19,0.85150670)
          (20,0.87148027)
          (21,0.88865520)
          (22,0.90341983)
          (23,0.91611562)
          (24,0.92703972)
          (25,0.93644870)
        };
        \addlegendentry{Board $A$}
        \addplot[red, thick, mark=square*, mark size=1.5pt] coordinates {
          (3,0.01502486)
          (4,0.05045260)
          (5,0.10538827)
          (6,0.17560372)
          (7,0.25565366)
          (8,0.34019693)
          (9,0.42469712)
          (10,0.50569674)
          (11,0.58082925)
          (12,0.64869046)
          (13,0.70865091)
          (14,0.76065903)
          (15,0.80506226)
          (16,0.84245830)
          (17,0.87357938)
          (18,0.89920748)
          (19,0.92011589)
          (20,0.93703184)
          (21,0.95061515)
          (22,0.96144840)
          (23,0.97003506)
          (24,0.97680266)
          (25,0.98210898)
        };
        \addlegendentry{Board $B$}
      \end{axis}
    \end{tikzpicture}
    \caption{Solo play: independent Markov chain}
    \label{fig:solo_cdf}
  \end{subfigure}\\
  \begin{subfigure}[b]{0.48\textwidth}
    \centering
    \begin{tikzpicture}
      \begin{axis}[
          width=\linewidth,
          title={PDF},
          xlabel={Roll number ($k$)},
          ylabel={$P(T=k)$},
          xmin=1, xmax=25,
          ymin=0, ymax=0.10,
          ytick={0,0.02,0.04,0.06,0.08,0.10},
          yticklabels={0,0.02,0.04,0.06,0.08,0.10},scaled y ticks=false,
          grid=both,
          major grid style={dotted, gray!50, line width=0.5pt},
          minor grid style={dotted, gray!20},
          minor tick num=1,
          legend style={at={(0.97,0.97)}, anchor=north east, font=\small},
        ]
        \addplot[blue, thick, mark=*, mark size=1.5pt] coordinates {
          (3,0.01680384) (4,0.03623685) (5,0.05115614) (6,0.05925428)
          (7,0.06095201) (8,0.05784595) (9,0.05175921) (10,0.04426466)
          (11,0.03651611) (12,0.02924945) (13,0.02285939) (14,0.01749583)
          (15,0.01315186) (16,0.00973262) (17,0.00710365) (18,0.00512174)
          (19,0.00365262) (20,0.00257939) (21,0.00180537) (22,0.00125343)
          (23,0.00086383) (24,0.00059130) (25,0.00040223)
        };
        \addlegendentry{Board $A$}
        \addplot[red, thick, mark=square*, mark size=1.5pt] coordinates {
          (3,0.01502486) (4,0.03285572) (5,0.04652899) (6,0.05373248)
          (7,0.05491045) (8,0.05166685) (9,0.04578490) (10,0.03875776)
          (11,0.03164310) (12,0.02508582) (13,0.01940775) (14,0.01470840)
          (15,0.01095168) (16,0.00803044) (17,0.00580984) (18,0.00415367)
          (19,0.00293837) (20,0.00205901) (21,0.00143051) (22,0.00098616)
          (23,0.00067504) (24,0.00045908) (25,0.00031036)
        };
        \addlegendentry{Board $B$}
      \end{axis}
    \end{tikzpicture}
    \caption{Head-to-Head play: joint Markov chain}
    \label{fig:joint}
  \end{subfigure}
  \begin{subfigure}[b]{0.48\textwidth}
    \centering
    \begin{tikzpicture}
      \begin{axis}[
          width=\linewidth,
          title={CDF},
          xlabel={Roll number ($k$)},
          ylabel={$P(T \le k)$},
          xmin=1, xmax=25,
          ytick={0,0.25,0.5,0.75,1},
          grid=both,
          major grid style={dotted, gray!50, line width=0.5pt},
          minor grid style={dotted, gray!20},
          minor tick num=1,
          legend style={at={(0.97,0.37)}, anchor=north   east, font=\small},
        ]
        \addplot[blue, thick, mark=*, mark size=1.5pt] coordinates {
          (3,0.01680384)
          (4,0.05304069)
          (5,0.10419683)
          (6,0.16345111)
          (7,0.22440312)
          (8,0.28224907)
          (9,0.33400828)
          (10,0.37827294)
          (11,0.41478905)
          (12,0.44403850)
          (13,0.46689789)
          (14,0.48439372)
          (15,0.49754558)
          (16,0.50727820)
          (17,0.51438185)
          (18,0.51950359)
          (19,0.52315621)
          (20,0.52573560)
          (21,0.52754097)
          (22,0.52879440)
          (23,0.52965823)
          (24,0.53024953)
          (25,0.53065176)
        };
        \addlegendentry{Board $A$}
        \addplot[red, thick, mark=square*, mark size=1.5pt] coordinates {
          (3,0.01502486)
          (4,0.04788058)
          (5,0.09440957)
          (6,0.14814205)
          (7,0.20305250)
          (8,0.25471935)
          (9,0.30050425)
          (10,0.33926201)
          (11,0.37090511)
          (12,0.39599093)
          (13,0.41539868)
          (14,0.43010708)
          (15,0.44105876)
          (16,0.44908920)
          (17,0.45489904)
          (18,0.45905271)
          (19,0.46199108)
          (20,0.46405009)
          (21,0.46548060)
          (22,0.46646676)
          (23,0.46714180)
          (24,0.46760088)
          (25,0.46791124)
        };
        \addlegendentry{Board $B$}
      \end{axis}
    \end{tikzpicture}
    \caption{Head-to-Head play: joint Markov chain}
    \label{fig:joint_cdf}
  \end{subfigure}
  \caption{(a)--(b) In solo play, modeled by independent Markov
    Chains, board~$B$ appears more likely to win. In fact, $P(T_B = k)
    > P(T_A = k)$ for $6 \le k \le 16$. (c)--(d) In a head-to-head
    matchup where both boards see the same dice rolls, modeled by
    joint Markov Chains, $P(T_A = k) > P(T_B = k)$ for $k\le 25$. In
    fact, $P(T_A < T_B) \approx 0.531$.}
  \label{fig:h2h}
\end{figure}
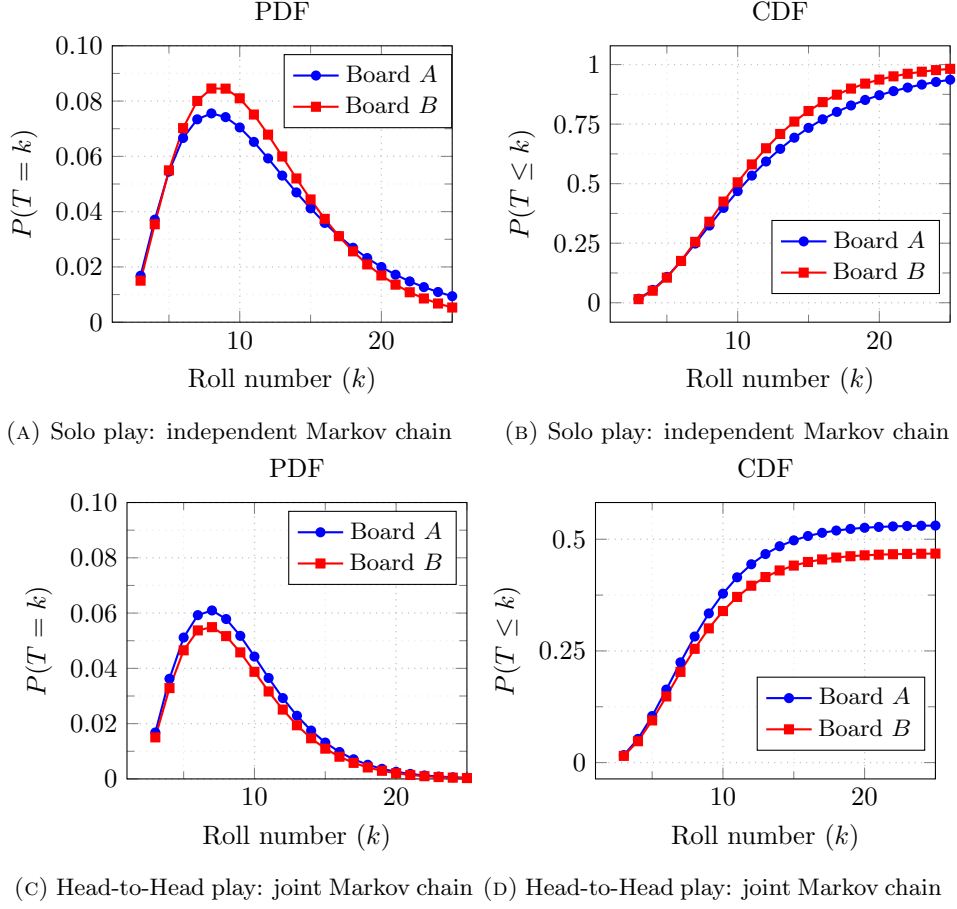

As a preliminary to constructing nontransitive triples of boards in
\cref{sec:nontransitive}, we observe that a board with a smaller
expected completion time is \emph{not} always favored in a head to head
matchup.  As a concrete example, consider boards $A$ and $B$ in
\cref{fig:h2h_boards}. Board $A$ is a Latin square of $\{6,7,9\}$, and
board $B$ uses only $\{6,7\}$.
\begin{figure}[htbp]
  \centering
  \begin{subfigure}[b]{0.4\textwidth}
    \centering
    \bingo{9}{6}{7}{7}{9}{6}{6}{7}{9} 
    \caption{Board $A$}
  \end{subfigure}
  \begin{subfigure}[b]{0.4\textwidth}
    \centering
    \bingo{6}{7}{6}{7}{7}{7}{6}{6}{6} 
    \caption{Board $B$}
  \end{subfigure}
  \caption{Two bingo boards considered for a head to head competition,
    with $V_A(\varnothing) \approx 12.77$ and $V_B(\varnothing)
    \approx 11.37$.}
  \label{fig:h2h_boards}
\end{figure}
In expectation, board $A$ requires nearly $1.4$ additional rolls to obtain a bingo in solo play, yet in a head-to-head,
board $A$ wins more than 53\% of the time,
\begin{align*}
  P(T_A < T_B) \approx 0.5315, \qquad P(T_B < T_A) \approx 0.4685,
  \qquad
  P(T_A = T_B) = 0.
\end{align*}
This surprising result can be explained by considering \cref{fig:h2h},
which shows the Probability Density Function (PDF) (probability for
winning on roll $k$) and the corresponding Cumulative Distribution
Function (CDF) for both the solo play (independent Markov chain) and
head-to-head play (joint Markov chain).
The contrast between the solo and head-to-head play is striking.  In
solo play, \cref{fig:solo,fig:solo_cdf}, board~$B$ dominates with
$P(T_B = k) > P(T_A = k)$ for $6 \le k \le 16$.  One would naturally
predict that board $B$ wins in a head-to-head competition. However, in
a head-to-head competition, \cref{fig:joint,fig:joint_cdf}, $P(T_A =
k) > P(T_B = k)$ for $k\le 25$.

This reversal occurs because the solo play involves independent Markov
chains, ignoring the correlation created by shared dice.  Dice sums of
$v=6$ or $v=7$ advance \emph{both} boards, roughly canceling out, but
a dice sum of $v=9$ (probability $4/36$) advances only board~$A$.
Board~$B$ has no exclusive sums -- every roll that helps $B$ also
helps~$A$.  This one-sided advantage manifests at almost  every roll number in
the joint game. 
This ``reversal'' makes it plausible that a trio of
nontransitive boards exists. A computational search identified
hundreds of examples; we present one such example in
\cref{sec:nontransitive}.

\section{Nontransitive boards}
\label{sec:nontransitive}

Recall from \cref{sec:h2h} that, for two boards $X$ and $Y$, we say that
\textbf{$X$ beats $Y$} if
\begin{align*}
  P(T_X<T_Y)>P(T_Y<T_X),
\end{align*}
where both players see the same sequence of dice rolls and each
follows their solo-optimal strategy.

It is useful to distinguish three related but different notions of
comparison: the solo expected completion time $V_X(\varnothing)$ of a
single board $X$; the head-to-head win probability when two boards see
the same dice rolls and each player follows their solo-optimal
strategy; and the fully strategic head-to-head game in which each
player may adapt choices to the opponent's board and current state.
In this section we use the second notion.  This comparison need not be
transitive.

We define a \textbf{nontransitive triple} to be three boards $C$,
$D$, and $E$ such that
\begin{align*}
  C \succ D,\qquad D \succ E,\qquad E \succ C,
\end{align*}
where $X\succ Y$ means that $X$ beats $Y$ in the above head-to-head
sense.  This is analogous to nontransitive dice, popularized by
Gardner \cite{Gardner1970NontransitiveDice}.  A smaller example of
nontransitive dice bingo boards appears in Gardner's book
\cite{Gardner1988TimeTravel}, where it is attributed to Knuth; see
also \cite{stackexchange_knuth_bingo}.  That example uses $2\times 2$
boards, a single die, and a winning condition based only on completing
a row.

We conducted a computational search over several hundred candidate
boards and found many nontransitive triples.  The example in
\cref{fig:nontransitive} uses only values from $\{5,6,7,8,9\}$:
\begin{figure}[htbp]
  \centering
  \begin{subfigure}[b]{0.3\textwidth}
    \centering
    \bingo{7}{7}{7}{6}{6}{6}{6}{7}{6}
    \caption{Board $C$}
  \end{subfigure}
  \begin{subfigure}[b]{0.3\textwidth}
    \centering
    \bingo{7}{5}{9}{9}{7}{5}{5}{9}{7}
    \caption{Board $D$}
  \end{subfigure}
  \begin{subfigure}[b]{0.3\textwidth}
    \centering
    \bingo{9}{7}{9}{9}{9}{9}{9}{6}{7}
    \caption{Board $E$}
  \end{subfigure}
  \caption{Three nontransitive boards with
    $V_C(\varnothing) \approx 11.31$,
    $V_D(\varnothing) \approx 11.34$, and
    $V_E(\varnothing) \approx 12.39$.}
  \label{fig:nontransitive}
\end{figure}
The head-to-head probabilities, computed in exact arithmetic as in
\cref{sec:h2h}, are shown in \cref{tbl:nontransitive}.
\begin{table}[h!]
  \centering
  \renewcommand{\arraystretch}{1.2}
  \begin{tabular}{lcccl}
    \hline
    \begin{tabular}[c]{@{}c@{}}
      Matchup\\
      $X$ vs.\ $Y$
    \end{tabular}
    & $P(X \text{ wins})$ & $P(Y\text{ wins})$ & $P(\text{tie})$ & Winner \\
    \hline
    $C$ vs.\ $D$ & 0.4281 & 0.3680 & 0.2039 & $C$ beats $D$ \\
    $D$ vs.\ $E$ & 0.3851 & 0.3177 & 0.2972 & $D$ beats $E$ \\
    $E$ vs.\ $C$ & 0.4948 & 0.4207 & 0.0845 & $E$ beats $C$ \\
    \hline
  \end{tabular}
  \caption{Head-to-head outcomes for the nontransitive triple in
    \cref{fig:nontransitive}.  In each row, board $X$ is favored over
    board $Y$ when $P(X\text{ wins})>P(Y\text{ wins})$, producing the
    cycle $C \succ D \succ E \succ C$.}
  \label{tbl:nontransitive}
\end{table}

Thus the head-to-head relation among these boards is genuinely cyclic:
\begin{align*}
  C \succ D \succ E \succ C.
\end{align*}
The margins are about $6$--$7\%$ in each matchup, after ignoring ties.
Notice also that board $E$ has the worst solo expected completion
time, $V_E(\varnothing)\approx 12.39$, while $V_C(\varnothing)\approx
11.31$; nevertheless, board $E$ beats board $C$ convincingly.  This
gives another instance of the phenomenon from \cref{sec:slower}: a
board that is slower in solo play can still be favored in a
head-to-head race.

\subsection{The role of variance}

One useful way to view the nontransitivity is through the spread of
the completion-time distributions. The standard deviations of the
three random variables $T_C,T_D, T_E$ are approximately
\begin{align*}
  \sigma_C \approx 5.32,\qquad
  \sigma_D \approx 5.63,\qquad
  \sigma_E \approx 7.52.
\end{align*}
Thus board $E$ has about $41\%$ larger standard deviation than board
$C$.  Its completion-time distribution is more spread out, with more
mass in both the left tail and the right tail.  In a head-to-head
race, the left tail is especially important, because unusually early
completions can decide the match before the opponent's better average
performance has time to appear.

This is the same broad phenomenon that appears in nontransitive dice:
when random variables are compared by ``who finishes first'' rather
than by expected value, the resulting ordering need not be transitive.
In dice bingo, however, the effect is amplified by the correlation
induced by the shared dice sequence.
\begin{remark}
  In our head-to-head matchups, \cref{sec:h2h,sec:nontransitive}, we
  assume each player follows their own optimal strategy, ignorant of
  the other boards. When the optimal strategy has multiple cells
  achieving the same minimum in \cref{eq:bellman}, the solo expected
  time is unaffected so they can choose at random from the equally
  fast options. However different choices could in principle yield
  different head-to-head probabilities.  For the boards presented
  here, we verified computationally that all tie-breaking rules
  produce identical head-to-head results. In general one might decide
  that the player should choose at random from among equally optimal
  choices, or simply mark the smallest numbered box from among equal
  choices.
\end{remark}

\begin{remark}
  Ties in the optimal strategy often arise from symmetries of the
  board, for example in board A if a six is rolled the optimal
  strategy is to mark the center. If a second six is rolled, either
  bottom corner is optimal and clearly that choice will not affect the
  winner in a head to head matchup since the two choices are perfectly
  symmetric. It is theoretically possible for an arithmetic
  coincidence to give two equally optimal choices that do not arise
  from any symmetry.
\end{remark}

\section{Game theory in head-to-head play}
\label{sec:nash}

\subsection{When strategy depends on the opponent}

In the preceding sections, both in our single board computations and
the head-to-head competitions for our nontransitive boards, each
player follows a strategy that minimizes their own expected completion
time, without regard to the opponent's board or progress.  In
head-to-head play, however, the objective is not to finish quickly but
to finish \emph{before the opponent}. These are different goals, and a
simple example shows just how dramatically they can diverge.

\subsection{A motivating example}

Consider the following position.  Player~2's board has cells 1 and~2
already marked, and is a single 7 away from completing the first row.
\begin{figure}[htbp]
  \centering
  \begin{subfigure}[b]{0.48\textwidth}
    \centering
    \bingo{7}{$\times$}{6}{9}{11}{8}{4}{6}{3} 
    \caption{Player 1: cell 1 marked}
  \end{subfigure} 
  \begin{subfigure}[b]{0.48\textwidth}
    \centering
    \bingo{7}{$\times$}{$\times$}{8}{9}{10}{11}{3}{5}
   \caption{Player~2: cells $\{1,2\}$ marked}
  \end{subfigure}
\end{figure}
Now suppose a dice sum $v=6$ is rolled.  Player~1 has two options of
cells to mark since $b_2=b_7=6$.
\begin{itemize}
\item \textbf{Cell~2}: The first row becomes two-thirds complete,
  needing one more~\textbf{7}---the same sum Player~2 needs.
\item \textbf{Cell~7}: The second column becomes two-thirds complete,
  needing one more~\textbf{11} -- a sum Player~2 does not need.
\end{itemize}

In solo play, the optimal choice is cell~2: the dice sum $v=7$ is
rolled with probability $6/36$, while the dice sum $v=11$ is rolled
with probability $2/36$.  The line through cell~2 completes three
times faster. But in the head-to-head game, cell~2 is a trap.  When
a~7 is eventually rolled, player~2 completes row~0 on \emph{their}
board -- and so does player~1.  The result is a tie.  The 7-line,
which is attractive in isolation, will never produce a \emph{win}
against this particular opponent. Thus, if our goal is to maximize the
chance of winning, player~1 should mark the~6 in cell~7. Then the game
becomes a race between the first occurrence of sum 11 and the first occurrence of sum 7, assuming no other winning
lines surface; player~1 needs a dice sum $v=11$ (rolled with
probability $2/36$) while player~2 needs a dice sum $v=7$ (rolled with
probability $6/36$).  Player~1's chance of winning this race is
\begin{align*}
  \frac{P(11)}{P(11) + P(7)} = \frac{2}{2+6} = \frac{1}{4}.
\end{align*}
A one-in-four chance of winning is far better than the zero percent
offered by the solo-optimal choice.

\begin{remark}
  This analysis focuses on the dominant line from each choice. In the
  full game, other winning lines may develop over subsequent rolls, so
  for example player 1 might end up winning with $6-8-3$ in the third
  column.  But the principle is clear: sharing winning sums with the
  opponent converts potential wins into ties, and a player who can see
  the opponent's board and wants to maximize their winning chances
  should sometimes choose a \emph{slower} line to avoid this.
\end{remark}

This observation opens the door to game theory: if both players adapt
their choices based on the opponent's board, then each player's
optimal decision may depend on the other's.  To analyze such
interactions, we use the concept of a \emph{Nash equilibrium}.

\begin{definition}
  In a two-player game where player~1 has strategies $A$ and $B$ and
  player~2 has strategies $X$ and $Y$, a \emph{Nash equilibrium} is a
  pair of strategies such that neither player can improve their payoff
  by unilaterally switching, with the other player's strategy held
  fixed.
\end{definition}
We can think of a Nash equilibrium as a stable set of strategies:
neither player has any incentive to change their strategy upon
learning the other's strategy. A Nash equilibrium can be \emph{pure}
(each player commits to a single choice) or \emph{mixed} (each player
randomizes over their choices with specified probabilities).  Nash's
theorem guarantees that every finite non-cooperative game has at least
one Nash equilibrium, though not necessarily a pure one. The classic
example of a game with no pure equilibrium is rock-paper-scissors: if
one player specifies a single choice, the other has a better
response. The unique Nash equilibrium is for both players to randomly
choose rock, paper or scissors with probability $1/3$. While this
strategy makes one unbeatable in the long run, the existence of
professional rock-paper-scissors leagues suggests psychology can also
play a role in this game!

\subsection{A position with no pure equilibrium}
We now exhibit a Dice Bingo position where the optimal head-to-head
play requires randomization.  Consider a head-to-head matchup between
the two boards in \cref{fig:nash_boards},
\begin{figure}[htbp]
  \centering
  \begin{subfigure}[b]{0.48\textwidth}
    \centering
    \bingo{7}{4}{6}{9}{7}{8}{6}{9}{5}
    \caption{Player 1 board}
  \end{subfigure}
  \begin{subfigure}[b]{0.48\textwidth}
    \centering
    \bingo{6}{9}{2}{3}{5}{5}{8}{10}{6} 
    \caption{Player 2 board}
  \end{subfigure}
  \caption{Two bingo boards where optimal head-to-head play requires
    randomization.}
  \label{fig:nash_boards}
\end{figure}
where the game begins with both boards unmarked.  Suppose the first
three shared rolls are $v=2$, $v=3$, and $v=10$.  Player~1's board
contains none of these values, so Player~1's board remains unmarked.
Player~2 marks three cells (see \cref{fig:nash_step}).
\begin{figure}[htbp]
  \centering
  \begin{subfigure}[b]{0.48\textwidth}
    \centering
    \bingo{7}{4}{6}{9}{7}{8}{6}{9}{5} 
  \end{subfigure} 
  \begin{subfigure}[b]{0.48\textwidth}
    \centering
    \bingo{6}{9}{$\times$}{$\times$}{5}{5}{8}{$\times$}{6} 
  \end{subfigure}
  \caption{After shared rolls of $2$, $3$ and $10$, player 2 has
    marked three cells, player 1's board has no marks.}
  \label{fig:nash_step}
\end{figure}
Suppose that the next roll is $v=6$.  Player~1 has two options of
cells to mark since $b_2 = b_6 = 6$. Player~2 also has two options of
cells to mark since $b_0 = b_8 = 6$.  The choices are simultaneous.

Denote player~1's strategies as
\begin{align*}
  A = \text{mark cell }2, \qquad B = \text{mark cell }6,
\end{align*}
and denote player~2's strategies as
\begin{align*}
  X = \text{mark cell }0, \qquad Y = \text{mark cell }8.
\end{align*}

After the simultaneous choice, the game continues with future shared
rolls, each player following their solo-optimal strategy from the
resulting state.  We compute the win probabilities in exact arithmetic
using the joint Markov chain discussed in \cref{sec:h2h} for each of
the four pairs of choices. The win probabilities are reported in a
{\em payoff matrix}, where each cell records $\bigl(P(\text{player 1
  wins}),\; P(\text{player 2 wins})\bigr)$, with the remainder being
the tie probability. The payoff matrix for our specific example is
given in \cref{tbl:nash_exact} in exact arithmetic; a numeric
approximation is given in \cref{tbl:nash_approx}.
\begin{table}[h!]
  \centering
  \renewcommand{\arraystretch}{2.0}
  \begin{tabular}{r|c|c|}
    \multicolumn{1}{c}{} & \multicolumn{1}{c}{$X$ (cell 0)} & \multicolumn{1}{c}{$Y$ (cell 8)} \\
    \cline{2-3}
    $A$ (cell 2) &
    $\left(\dfrac{17279}{86184},\;\dfrac{10186789}{13961808}\right)$ &
    $\left(\dfrac{656700799}{3411916830},\;\dfrac{3462999829}{4616122770}\right)$ \\[4pt]
    \cline{2-3}
    $B$ (cell 6) &
    $\left(\dfrac{155}{912},\;\dfrac{110569}{147744}\right)$ &
    $\left(\dfrac{22214689}{114686280},\;\dfrac{281297}{387828}\right)$ \\[4pt]
    \cline{2-3}
  \end{tabular}
  \caption{Payoff matrix for the four pairs of choices.  Each cell
    records $\bigl(P(\text{player 1 wins}),\; P(\text{player 2
      wins})\bigr)$, with the remainder being the tie
    probability. Computations are done in exact arithmetic.}
  \label{tbl:nash_exact}
\end{table}
\begin{table}[h!]
  \centering
  \renewcommand{\arraystretch}{1.5}
  \begin{tabular}{r|c|c|}
    \multicolumn{1}{c}{} & \multicolumn{1}{c}{$X$ (cell 0)} & \multicolumn{1}{c}{$Y$ (cell 8)} \\
    \cline{2-3}
    $A$ (cell 2) & $(0.2005,\; 0.7296)$ & $(0.1925,\; 0.7502)$ \\
    \cline{2-3}
    $B$ (cell 6) & $(0.1700,\; 0.7484)$ & $(0.1937,\; 0.7253)$ \\
    \cline{2-3}
  \end{tabular}
  \caption{A floating point approximation to \cref{tbl:nash_exact}.}
  \label{tbl:nash_approx}
\end{table}

Player~2 is the heavy favorite regardless of the choices (winning
roughly $73$--$75\%$ of the time), but the precise probabilities
depend on both players' decisions.  We proceed to check each player's
best response to each of their opponent's strategies.

\begin{itemize}
\item
  \textbf{Player~1's best responses:} If Player~2 chooses $X$,
  Player~1 compares the first coordinates in the $X$-column: $0.2005 >
  0.1700$, so Player~1 prefers $A$.  If Player~2 chooses $Y$, Player~1
  compares the $Y$-column: $0.1937 > 0.1925$, so Player~1 prefers $B$.

\item
  \textbf{Player~2's best responses:} If Player~1 chooses $A$,
  Player~2 compares the second coordinates in the $A$-row: $0.7502 >
  0.7296$, so Player~2 prefers $Y$.  If Player~1 chooses $B$, Player~2
  compares the $B$-row: $0.7484 > 0.7253$, so Player~2 prefers $X$.
\end{itemize}
Thus no cell of the payoff matrix is a mutual best response and thus
this position has \textbf{no pure Nash equilibrium}; any fixed choice
of strategies (cells in the grid) leaves one player wanting to make a
change.

\subsection{The mixed Nash equilibrium}

Since no pure equilibrium exists, we instead compute the {\em mixed
  Nash equilibrium}. We assign probabilities for each player's two
choices and then solve the linear equations imposed by the condition
that neither player can improve their strategy after seeing the other
player's strategy; we say each player is \textit{indifferent} to their
opponent's strategy. So let
\begin{align*}
  p = P(\text{Player~1 chooses } A), \qquad
  q = P(\text{Player~2 chooses } X).
\end{align*}
At equilibrium, each player must be indifferent between their two
strategies.

\textbf{Player~1's indifference condition.}  Player~1's winning
probability under $A$ is
\begin{align*}    
  W_1(A) = q \cdot \frac{17279}{86184} + (1-q) \cdot
  \frac{656700799}{3411916830}\,,
\end{align*}
and under $B$:
\begin{align*}
  W_1(B) = q \cdot \frac{155}{912} + (1-q) \cdot \frac{22214689}{114686280}\,.
\end{align*}
Setting $W_1(A) = W_1(B)$ and solving gives
\begin{align*}
  q = \frac{352522}{9125389} \approx 0.0386.
\end{align*}

\textbf{Player~2's indifference condition.}  Player~2's winning
probability under $X$ is
\begin{align*}
  W_2(X) = p \cdot \frac{10186789}{13961808} +
  (1-p) \cdot \frac{110569}{147744}\,,
\end{align*}
and under $Y$:
\begin{align*}
  W_2(Y) = p \cdot \frac{3462999829}{4616122770} +
  (1-p) \cdot \frac{281297}{387828}\,.
\end{align*}
Setting $W_2(X) = W_2(Y)$ and solving gives
\begin{align*}
  p = \frac{29891145}{56555767} \approx 0.5285.
\end{align*}

The mixed Nash equilibrium is therefore: Player~1 marks cell~2 with
probability $p\approx 53\%$ and cell~6 with probability $1-p\approx
47\%$, while Player~2 marks cell~8 with probability $1-q\approx 96\%$
and cell~0 with probability $q\approx 4\%$.  The strong asymmetry in
Player~2's mixture reflects the fact that cell~8 is nearly always the
right choice; Player~2 mixes in cell~0 just often enough to keep
Player~1 from exploiting the predictability. Either player can
announce their mixed strategy and the other player can do no better
than their own mixed strategy.

\begin{remark}
  This example demonstrates that Dice Bingo is not merely a
  single-player optimization problem.  Once players account for the
  opponent's board, a single shared roll can create a genuine
  strategic subgame requiring game-theoretic reasoning.  The practical
  takeaway: in head-to-head play, occasionally deviating from
  solo-optimal play -- in a deliberately unpredictable way -- can be
  part of the optimal strategy. It also means trying to calculate win
  probabilities or to find nontransitive boards using the version of
  head-to-head matchups where opponents utilize game theory is a much
  more complex problem. We would need to analyze the Nash equilibrium
  choices at each state $(S_A,S_B)$.
\end{remark}

\section{Conclusions and future exploration}

We have presented Dice Bingo as a rich setting for exploring concepts
from linear algebra, probability, finite Markov chains, Bellman
equations, group actions and game theory. We model optimal play for a
board as a finite Markov decision process, solving the associated
Bellman equations to identify a unique optimal board up to natural
symmetries.  We then examine head-to-head competition and show that a
board with a worse expected completion time can be favored in
head-to-head competition. Indeed, we discover nontransitive sets of
bingo boards, where board $A$ is favored against board $B$, board $B$
is favored against board $C$, and board $C$ is favored against board
$A$. Finally, we study strategic play where players adapt their choice
to their opponent's board. In this setting, we exhibit a position with
no pure Nash equilibrium and compute an explicit mixed Nash
equilibrium. One could also replace the head-to-head with $n$ players
competing to obtain the first bingo.

A GitHub repository containing Python scripts related to this
manuscript is available at
\url{https://github.com/ongbw/dice-bingo}. If preferred, an
interactive Python notebook can be launched to replicate results in
this manuscript by visiting
\url{https://mybinder.org/v2/gh/ongbw/dice-bingo/main?filepath=dice-bingo.ipynb}.
To experiment with dice bingo, readers can visit this interactive html webpage,
\url{https://ongbw.github.io/dice-bingo/dice-bingo.html}.

There are numerous intriguing variations that remain
unexplored. Beyond using larger bingo boards, one could investigate
boards of different shapes with unconventional winning patterns, as
well as experimenting with different probabilities for the bingo cell
entries—such as sums of more than two dice, sums from other
loaded dice, or alternative formats like product bingo
\cite{Jackson2013Bingo}.

\section*{Acknowledgements }
The High-Performance Computing Shared Facility (Superior and/or
Portage) at Michigan Technological University was used in obtaining
results presented in this manuscript.

\bibliographystyle{vancouver}   
\bibliography{references}   
\end{document}